\documentclass[letterpaper11pt]{article}

\usepackage[utf8]{inputenc}
\usepackage[T1]{fontenc}
\usepackage[breaklinks]{hyperref}
\usepackage[margin=1in]{geometry}
\usepackage{amssymb}
\usepackage{amsmath}
\usepackage{dsfont}
\usepackage{mathtools}
\usepackage{authblk}
\usepackage{xcolor}
\usepackage{dsfont}
\usepackage{subcaption}
\usepackage{wrapfig}
\usepackage{amsthm}
\usepackage{xurl}

\makeatletter
\providecommand\mainlastpage{0}
\providecommand\applastpage{0}
\AtEndDocument{%
  \immediate\write\@auxout{%
    \string\csxdef{lastpage}{\mainlastpage}%
    \string\csxdef{applastpage}{\thepage}}}
\makeatother

\newcommand{\E}{\mathbb{E}}

\newcommand{\var}{\text{Var}}

\newcommand{\cov}{\text{Cov}}
\newcommand{\corr}{\text{Corr}}

\newtheorem{theorem}{Theorem}

\newtheorem{defin}{Definition}

\usepackage[authoryear,longnamesfirst]{natbib}

\title{\LARGE \textbf{Control Variates to reduce the variance of the empirical Moment Generating Function: an extreme value perspective}}
\author[1,2]{Louison Bocquet-Nouaille \thanks{Corresponding author \\ Emails: louison.bocquet\_nouaille@onera.fr, jerome.morio@onera.fr, benjamin.bobbia@isae-supaero.fr, philippe.naveau@lsce.ipsl.fr}}
\author[1,2]{Jérôme Morio}
\author[2]{Benjamin Bobbia}
\author[3]{Philippe Naveau}

\affil[1]{ONERA/DTIS, Université de Toulouse, F-31055 Toulouse}
\affil[2]{Fédération ENAC ISAE-SUPAERO ONERA, Université de Toulouse, 31000 Toulouse}
\affil[3]{LSCE, CEA, CNRS, UVSQ, Université Paris-Saclay, Gif-sur-Yvette, France}

\date{\vspace{-10pt} \small{Preprint - November 11, 2026}}

\begin{document}

\maketitle

\vspace{10pt}
Accurate estimation of Moment Generating Functions is fundamental in many applications, but empirical estimators have high variances. Within a control variate setup, Extreme Value Theory is leveraged to identify how extreme values affect the variance reduction method's optimality.
\vspace{20pt}

\section{Introduction}\label{sec:intro}
Moment Generating Functions (MGF) play a fundamental role across applied probability and statistics. For instance, the MGF can be used for making decisions under uncertainty \citep{collender1986alternative} as an alternative to mean-variance analysis, to construct statistical tests \citep{cabana2005using, meintanis2007kolmogorov}, or in importance sampling as exponential tilting \citep{fuh2024efficient}. 
Its mathematical definition is simple. Given any real valued random variable $T$, its MGF is defined as 
\begin{align*}
    M(u) = \E\left[\exp\left(u T\right)\right] ,
    \label{eq:mgf}
\end{align*}
provided this expectation exists for all $u \in (-h, h)$, $h > 0$. 
The MGF $M(u)$ uniquely determines the distribution of $T$.

In non-parametric settings, the MGF is naturally estimated with its empirical counterpart \citep{Gbur01011989}. More precisely, let ${\bf T}_n=(T_1, \dots, T_n)^\perp$ be $n$ i.i.d. copies of the random variable $T$. Then, the empirical estimate of $M(u)$, say EMGF, can be obtained as the following mean
\begin{align*}
    \widehat{M}_n(u) = \frac{1}{n}\sum_{i=1}^n \tilde{T}(u),
\end{align*}
where the random variable $\tilde{T}(u) = \exp\left(u T_i\right)$ denotes the exponentiated version of each $T_i$.

Exponentiating observations can lead to severe variance, especially for large values of $u$ or small sample sizes $n$. 
To mitigate this issue, several non-parametric improvements have been proposed, such as Bernstein polynomial smoothing \citep{Susam15092026} and kernel-type estimators \citep{BOUZEBDA2022109422}. However, these methods trade variance reduction for bias and require delicate tuning of kernels and parameters (the step-size and bandwidth in \citet{BOUZEBDA2022109422} or the smoothing degree in \citet{Susam15092026}). In contrast our approach preserves unbiasedness and handles general distributions, unlike Bernstein-based smoothing which is confined to compact supports.

In this work, we tailor a control variate approach \citep[Section 2, Chapter 5]{asmussen_variance-reduction_2007} to reduce the variance of $\widehat{M}_n(u)$, while keeping an unbiased estimate.

One fundamental hypothesis in this control variate framework is to assume that two samples with two different distributions, but strongly correlated, are at the disposal of the statistician. 
The core idea in this setup is to improve statistical inference of quantities of interest, like $\widehat{M}_n(u)$ here, based on the called {\it target} sample ${\bf T}_n$, by leveraging information from the second related {\it source} sample, say ${\bf S}_n=(S_1,\dots, S_n)^\perp$. This can be viewed as a form of transfer learning from the source ${\bf S}_n$ to the target ${\bf T}_n$ \citep{zhu_recent_2025}. 
Specifically, we assume to have access to $n$ i.i.d. copies from the bivariate random vector $\left(T,S\right)$. We further assume either that the expectation $G(u)=\E[\exp(u S)]$ is known, as in the exact control variate approach (see Definition~\ref{def:cv}), or that an additional sample of $m$ source observations $(S_{n+1},\dots, S_{n+m})$ is available to estimate $G(u)$, as in the approximate control variate approach (see Definition~\ref{def:acv}). Note that this framework naturally extends to multi-source settings \citep[pages 141-142]{asmussen_variance-reduction_2007}.

A novelty of this work is the study of the extremal behaviour of the proposed variance-reduced EMGF for large $u$, in particular identifying under which mathematical conditions a control variance setup can (or cannot) be effective at reducing the variance of estimators of $G(u)$. This result relies on adequately coupling Extreme Value Theory (EVT) concepts and control variate techniques. 

The paper is organized as follows. Section~\ref{sec:cv} presents the control variates variance-reduction method and its application to the EMGF. Section~\ref{sec:asymp} studies the asymptotic properties of the proposed estimator through EVT. Section~\ref{sec:simu} illustrates on simulated data the limits discussed in Section~\ref{sec:asymp}. For readability, all proofs are deferred to supplementary materials.

\section{Variance reduction with control variates}\label{sec:cv}
The control variates method \citep[Chapter 5, Section 2]{asmussen_variance-reduction_2007} adjusts the empirical mean estimator on the target data by the mean estimation error on the source data, weighted by a coefficient chosen to minimize the variance. By construction, this adjustment guarantees reduced variance whenever the target and source observations are linearly correlated, providing a powerful tool to improve MGF estimation if coupled samples are available.
The following definition clarifies this point.

\begin{defin}[Control variates MGF estimator] \label{def:cv}
Suppose that $n$ i.i.d. copies from the bivariate real-valued random vector $\left(T,S\right)$ are available.
If the source MGF 
$G(u)=\E\left[\exp\left(u S\right)\right]$ is known, but the target MGF $M(u)=\E\left[\exp\left(u T\right)\right]$ needs to be estimated, then the so-called control variate (CV) estimator can be defined as follows 
\begin{align*}
\widehat{M}^{CV}_n(u) = \widehat{M}_n(u) + w_{opt}(u) \left(G(u) - \widehat{G}_n(u)\right) ,
\end{align*}
where $\widehat{M}_n(u)$ and $\widehat{G}_n(u)$ represent the empirical means based on each exponentiated sample, respectively $\tilde{\bf T}_n(u)$ and $\tilde{\bf S}_n(u)$, with the convention $\tilde{x}=\exp\left(u x\right)$. The coefficient $w_{opt}(u)$ is chosen to be optimal in the sense that 
\begin{align*}
w_{opt}(u) := \underset{w \in \mathbb{R}}{\text{argmin}} \; \var\left(\widehat{M}_n(u) + w \left(G(u) - \widehat{G}_n(u)\right)\right) = \frac{\cov\left(\tilde{T}(u), \tilde{S}(u)\right)}{\var\left(\tilde{S}(u)\right)}.
\end{align*}
\end{defin}

By construction, $\widehat{M}^{CV}_n(u)$ is an unbiased estimator of $M(u)$. In addition, its variance can be viewed as a function of the baseline variance, i.e. $\var\left(\widehat{M}_n(u)\right)$, and the squared Pearson correlation via 
\begin{equation} \label{eq:var_red}
\var\left(\widehat{M}^{CV}_n(u)\right) = \left(1 - \corr\left(\tilde{T}(u),\tilde{S}(u)\right)^2\right) \var\left(\widehat{M}_n(u)\right).
\end{equation}
Consequently, the variance of $\widehat{M}^{CV}_n(u)$ never increases and its variance reduction is entirely governed by the magnitude of $\big|\corr\left(\tilde{T}(u),\tilde{S}(u)\right)\big|$. In the absence of correlation between the source and target samples, ${\bf S}_n$ and $ {\bf T}_n$, the control variate estimator reduces to the baseline, i.e. the EMGF defined by $\widehat{M}_n(u)$. Conversely, a variance gain is always observed when the source and target samples are correlated. 
In other words, the efficiency of the control variate method is governed by the factor $1 - \corr\left(\tilde{T}(u), \tilde{S}(u)\right)^2$, and insuring that extreme fluctuations in $\tilde{S}$ scale proportionally with those in $\tilde{T}$ for all values will be key when choosing the appropriate source control variate $\tilde{S}$. 
The strictly monotonic exponential function $\tilde{x}=\exp\left(u x\right)$, for any fixed non-null $u$, will play an important role as a strong correlation between $T$ and $S$ may or may not translate into high correlations between $\tilde{T}(u)$ and $\tilde{S}(u)$ for different $u$.

From a practical perspective, two hurdles quickly appear when computing $\widehat{M}^{CV}_n(u)$. A first limitation is to get the correct value of $w_{opt}(u)$. One solution is to arbitrarily choose a weight $w$ and plug it in our CV estimate to get $\widehat{M}^{CV}_n(u)$. This can hinder the variance reduction as the reduction occurs if and only if $w \in \left(\min(0, 2w_{opt}(u)), \max(0, 2w_{opt}(u))\right)$ \citep[Eq. (3)]{lavenberg_perspective_1981}. 
If $\widehat{w}_{opt}$ can be estimated from the data at hand, it has to be sufficiently concentrated around $w_{opt}(u)$ to insure to be in the interval $\left(\min(0, 2w_{opt}(u)), \max(0, 2w_{opt}(u))\right)$. Same-sample estimation of $w_{opt}(u)$ may also induce a plug-in bias \citep{kennedy_semiparametric_2023}, typically negligible in practice \citep{owen_variance_2013}.

The second hurdle to compute $\widehat{M}^{CV}_n(u)$ comes from the requirement of a known control variate MGF $G(u)$ in Definition \ref{def:cv}. This appears unrealistic in most applications.
To address this issue, the literature on control variates often makes an additional assumption concerning data availability. 
The so-called approximate control variates method \citep{gorodetsky_generalized_2020} assumes that $m$ additional source observations, say $(S_{n+1}, \dots, S_{n+m})^\perp$, are available to estimate $G(u)$.
The reasoning is justified whenever source data are easy to get, while target data are sparse \citep[see also][]{gorodetsky_generalized_2020}. 
Under this framework, a new estimator for $M(u)$ can be defined as described by the following definition.

\begin{defin}[Approximate control variates MGF estimator] \label{def:acv}
Under the same conditions as Definition \ref{def:cv}, but with an unknown MGF source $M(u)$ and an additional source sample of size $m$, say $(S_{n+1}, \dots, S_{n+m})^\perp$, one can define the approximate control variates estimator $\widehat{M}^{ACV}_n(u)$ of $M(u)$ as 
    \begin{align*}
    \widehat{M}^{ACV}_n(u) = \widehat{M}_n(u) + w_{opt}(u) \left(\widehat{G}_{n+m}(u) - \widehat{G}_n(u)\right) .
    \end{align*}
    The optimal coefficient $w_{opt}(u)$ minimizing the variance is the same as the one of $\widehat{M}^{CV}_n$ in Definition \ref{def:cv}.
\end{defin}
The estimate $\widehat{M}^{ACV}_n$ remains unbiased and its variance becomes 
\begin{eqnarray*}
    \var\left(\widehat{M}^{ACV}_n(u)\right)  &=& \var\left(\widehat{M}_n(u)\right) + \frac{m}{n+m} \left(\var\left(\widehat{M}^{CV}_n(u)\right)- \var\left(\widehat{M}_n(u)\right)\right),\\
    &=& 
    \left(1 - \frac{m}{n+m}\corr\left(\tilde{T}(u),\tilde{S}(u)\right)^2\right) \var\left(\widehat{M}_n(u)\right).
\end{eqnarray*}
Again, the correlation $\corr\left(\tilde{T}(u),\tilde{S}(u)\right)$ plays a fundamental role to assess the level of variance reduction. 
In this work, we investigate how this quantity behaves for large $u$. 
Similarly, one may wonder how the optimal coefficient $w_{opt}(u)$ changes when $u$ gets large as it governs the entire behaviour of the method. We concentrate our analysis on large values of $u$ because the small-$u$ regime is well-behaved, with control variates yielding standard variance reduction if $\tilde{T}$ and $\tilde{S}$ are correlated and if the employed coefficient is within the variance-reduction bounds $\left(\min(0, 2w_{opt}(u)), \max(0, 2w_{opt}(u))\right)$.
To address these problems, we rely on the field of EVT that is tailored to model tail behaviours \cite[see, e.g.][]{de_haan_extreme_2006}.

\section{Extremal behaviour of MGF control variate estimates}\label{sec:asymp}
Estimating the MGF functions $M(u)$ and $G(u)$ implicitly means that they must exist and consequently, all moments of $T$ and $S$ have to be finite. The existence of MGFs also provides important information about their upper tail behaviours. For example, neither $T$ nor $S$ can be Pareto, Cauchy or Student distributed as these distributions are heavy tailed, in particular they possess some infinite moments. To avoid such cases, we need to impose some conditions about the tail behaviours of $T$ and $S$. First, we assume that $T$ and $S$ have finite upper points, say $t^* = \sup\{x : F^T(x) < 1\}$ and $s^* = \sup\{x : F^S(x) < 1\}$ where $F^T(x)=P(T\leq x)$ and $F^S(x)=P(S\leq x)$, respectively. 
In addition, the speed of convergence towards these end points needs to be clearly quantified. Following classical EVT arguments, we also assume that 
\begin{align}\label{eq: RV condition}
\mathbb{P}\left(\frac{1}{t^*-T} > u\right) \underset{u \to +\infty}{\sim} u^{1 / \gamma^T} L_T(u) \quad \text{ and } \quad \mathbb{P}\left(\frac{1}{s^*-S} > u\right) \underset{u \to +\infty}{\sim} u^{1 / \gamma^S} L_S(u),
\end{align}
where $L_T(.)$ and $L_S(.)$ are both slowly varying function for large $u$, and the negative indices $\gamma^T$ and $\gamma^S$ drive the tail behaviours.
In EVT, such indices are called Extreme Value Index (EVI) and their negative nature make $T$ and $S$ members of the so-called Weibull max-stable domain of attraction. 
The last ingredient needed before studying the variance reduction of 
$\widehat{M}^{ACV}_n(u)$ for large $u$ is the description of the joint extremal dependence between $T$ and $S$. 

Precisely, we assume that 
\begin{align}\label{eq: MRV condition}
\underset{u \to +\infty}{lim} \frac{1}{g(u)}\mathbb{P}\left( \left(\frac{1}{u \left(t^*-T\right)},\frac{1}{u \left(s^*-S\right)}\right) \in B \right)= \nu(B), \ \text{ with } B = \left\{ (x,y)\in(0,\infty)^2; \frac{1}{x} + \frac{1}{y} < 1 \right\},
\end{align} 
where $g$ represents a positive function going to 0 at infinity, and $\nu(\cdot)$ a measure defined on the unit sphere. 
This is called a multivariate regular variation condition \citep{basrak2025foundations}.

\subsection{Extremal correlation}
As shown in \eqref{eq:var_red}, the theoretical variance reduction achieved by control variates ultimately depends directly on $1 - \corr\left(\tilde{T}(u), \tilde{S}(u)\right)^2$.
Understanding the extremal behaviour of this correlation will reveal which parameters influence the control variate adjustment and keep it effective.

\begin{theorem} \label{th:lim_corr}
Under Definitions \ref{def:cv} and \ref{def:acv}, conditions \eqref{eq: RV condition} and \eqref{eq: MRV condition}, and
if $\gamma^T = \gamma^S$ and $\nu\left( (0, \infty)^2\right)>0$, then the variance reduction factor $1 - \corr\left(\tilde{T}(u), \tilde{S}(u)\right)^2$ satisfies
\begin{equation} \label{eq:lim_corr}
    \underset{u \to +\infty}{lim} \left( 1 - \corr\left(\tilde{T}(u), \tilde{S}(u)\right)^2 \right) = 1 - \frac{\nu(B)^2}{c_T c_S} 2^{-2 / \gamma^T},
\end{equation}
where $c_T$ and $c_S$ are two positive constants. 
If either $\gamma^T \neq \gamma^S$ or $\nu\left( (0, \infty)^2\right)=0$, the above limit of the variance reduction factor becomes 1, meaning no variance reduction happens. 
\end{theorem}

The proof is detailed in supplementary material (Section~\ref{sec:proof_lim_corr}).
Theorem~\ref{th:lim_corr} establishes a fundamental theoretical boundary for statistical transfer in the proposed $\widehat{M}^{CV}_n(u)$ and $\widehat{M}^{ACV}_n(u)$ estimators: variance reduction under extreme parameters $u$ is asymptotically achievable if there is asymptotic dependence and if the target and source variables share identical EVIs. Indeed, if $\gamma^T \neq \gamma^S$, one tail systematically dominates on an exponential scale, and as $u \to +\infty$, the probability mass of the dominant variable concentrates at a rate that completely decorrelates $\tilde{T}(u)$ and $\tilde{S}(u)$. In this regime, control variates provide no asymptotic variance reduction over standard empirical mean estimation. However, when $\gamma^T = \gamma^S$, the joint extreme behaviour is preserved under exponential transformation and the limiting correlation converges to a non-zero constant, guaranteeing that the proposed $\widehat{M}^{CV}_n(u)$ and $\widehat{M}^{ACV}_n(u)$ estimators retain their variance-reduction property.

\subsection{Extremal optimal coefficient}
Building on Theorem~\ref{th:lim_corr} which established that asymptotic variance reduction is achievable if $\gamma^T = \gamma^S$ and $\nu\left( (0, \infty)^2\right) > 0$, we now focus our analysis of the optimal coefficient $w_{opt}(u)$ strictly on this viable transfer regime.

\begin{theorem}\label{th:lim_coef} Under the same conditions than Theorem \ref{th:lim_corr}, 
if $\gamma^T = \gamma^S$, $t^*=s^*$ and $\nu\left( (0, \infty)^2\right)>0$, then the optimal coefficient $w_{opt}(u)$ satisfies, for some positive constant $c_S$, 
\begin{equation}\label{eq:lim_coef}
    \underset{u \to +\infty}{lim} w_{opt}(u) = \frac{\nu(B)}{c_S}2^{-1 / \gamma^S},
\end{equation}
and otherwise, we have 
    \begin{align*}
        \underset{u \to +\infty}{lim} w_{opt}(u) = \begin{cases}
            + \infty & \text{if } t^* > s^*; \\
            0 & \text{if } t^* < s^*. \\
        \end{cases}
    \end{align*}
\end{theorem}

The proof is detailed in supplementary material (Section~\ref{sec:proof_lim_coef}). Alternative scenarios, such as different EVIs ($\gamma^T \neq \gamma^S$), or asymptotic independence ($\nu\left( (0, \infty)^2\right) = 0$), yield degenerate limits ($0$ or $\pm \infty$) and are comprehensively detailed along with their proofs in the supplementary material (Section~\ref{sec:proof_lim_coef}).

Under the condition $\gamma^T = \gamma^S$ and $\nu\left( (0, \infty)^2\right) > 0$, the limiting optimal coefficient $w_{opt}(u)$ converges to a stable finite non-zero constant if the endpoints match. This stability is a key prerequisite for practical variance reduction, ensuring that the optimal coefficient remains well-behaved in the extreme regime.

Ultimately, this limiting behaviour highlights the practical limits of control variates at high parameter values $u$ as the theoretical optimal coefficient $w_{opt}(u)$ is unknown and must be replaced by an estimate $\widehat{w}_{opt}$. A limit $w_{opt}(u) \to 0$ reflects no transfer, while $w_{opt}(u) \to +\infty$ reflects an unstable rescaling of mismatched endpoints. From an estimation standpoint, when $w_{opt}(u)$ vanishes or diverges, sample-based plug-in estimators $\widehat{w}_{opt}$ exhibit extreme variance as they rely on fourth-order moments that grow exponentially with $u$. This instability frequently pushes $\widehat{w}_{opt}$ outside the variance-reduction boundary $(0, 2w_{opt}(u))$ (or $(2w_{opt}(u), 0)$ if $w_{opt}(u) < 0$), causing the control variate method to inadvertently inflate, rather than reduce, the target variance.

\section{Simulation study}\label{sec:simu}
To illustrate the theoretical results of Theorems~\ref{th:lim_corr} and \ref{th:lim_coef}, the variance reduction factor $1 - \corr\left(\tilde{T}(u), \tilde{S}(u)\right)^2$ and the optimal coefficient $w_{opt}(u)$ are estimated using a sample of size $n = 1,000,000$. The samples follow Generalized Pareto marginals with scale parameter $\sigma=0.1$ and EVIs $(\gamma^T, \gamma^S) \in \{-0.5;-3.5\}^2$. Upper endpoints are set to $t^* \in \{0.5, 1, 1.5\}$ and $s^* = 1$. Asymptotic dependence is modelled by a Gumbel copula with parameter $\theta=3$. Extreme and contrasting EVI values, alongside a low dependence parameter, are chosen to accelerate convergence and avoid the computational burden of larger sample sizes. The parameter $u$ ranges from 1 to 150.

\begin{figure}[!h]
    \centering
    \includegraphics[width=\linewidth]{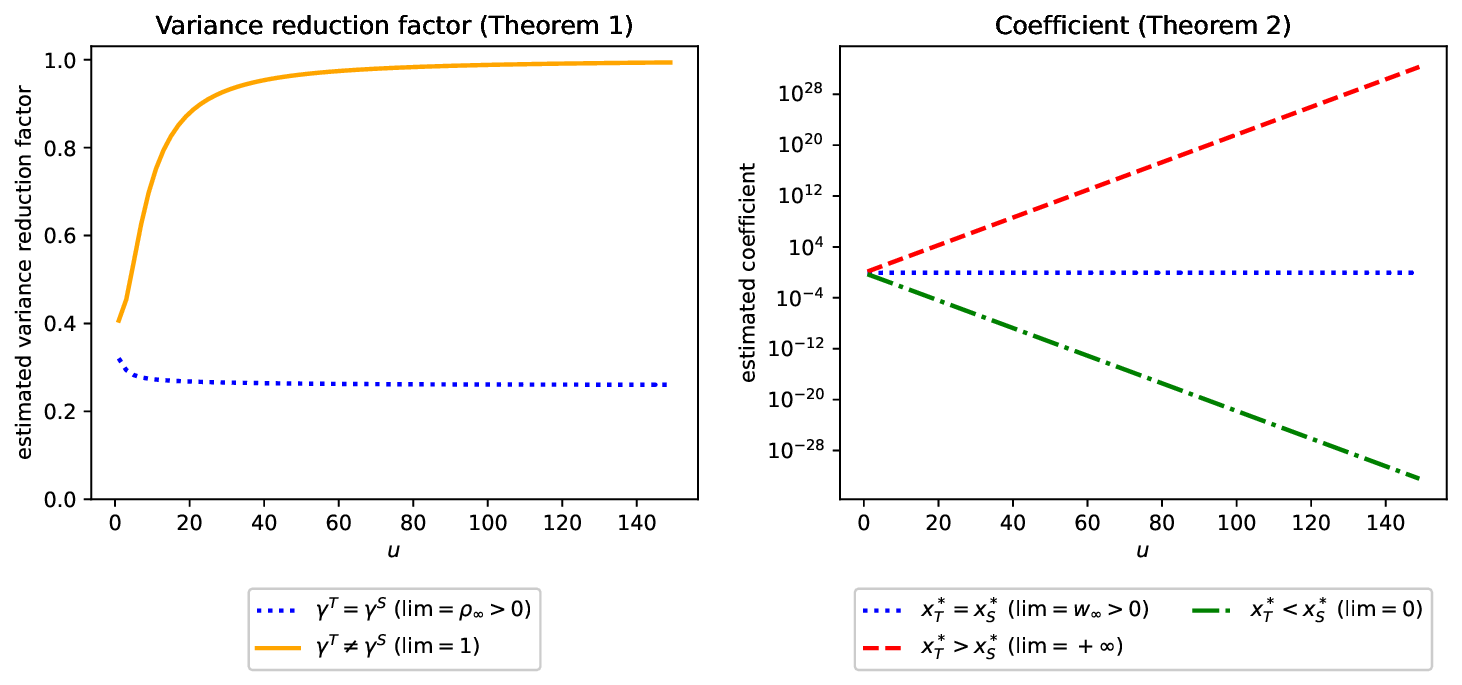}
    \caption{\textbf{Estimated variance reduction factor $1 - \corr\left(\tilde{T}(u), \tilde{S}(u)\right)^2$ (left) and estimated optimal coefficient $w_{opt}(u)$ (right) as functions of $u$.} For clarity, only one representative configuration per limiting regime is displayed. $\rho_{\infty}$ and $w_{\infty}$ denote the finite non-zero limits given in Eq.~\eqref{eq:lim_corr} (Theorem~\ref{th:lim_corr}) and Eq.~\eqref{eq:lim_coef} (Theorem~\ref{th:lim_coef}) respectively.}
    \label{fig:lim}
\end{figure}

Figure~\ref{fig:lim} displays the simulation results, capturing the limiting behaviours of both the correlation and the optimal coefficient as $u$ increases. All combinations of endpoints and EVIs are simulated, but only representative configurations are shown; the others exhibit the same limiting behaviour, as predicted by Theorems~\ref{th:lim_corr} and \ref{th:lim_coef}. The variance reduction factor $1 - \corr\left(\tilde{T}(u), \tilde{S}(u)\right)^2$ (left) is shown for equal endpoints only, since its limit does not depend on them, and the case $\gamma^T \neq \gamma^S$ is illustrated with $\gamma^T > \gamma^S$. For the coefficient (right), endpoint mismatches are illustrated with $\gamma^T = \gamma^S$.

\section{Discussion}\label{sec:conclu}
To conclude, if correlated paired samples are available, control variates offer a powerful approach to reduce the variance of the EMGF without introducing bias. However, our asymptotic analysis using EVT establishes critical boundaries when the parameter $u$ becomes large. Specifically, even with knowledge of the optimal coefficient, no variance reduction can be obtained for large $u$ if the target and source EVIs do not match. Furthermore, our study of the optimal control variate coefficient's limit reveals severe instabilities when endpoints are misaligned, hindering the coefficient's estimation which can lead to a variance increase.

We analysed the asymptotic properties as $u \to +\infty$ but the limiting behaviour for $u \to -\infty$ follows analogously by considering the negated variables $T' = -T$ and $S' = -S$, where left endpoints $u'_* = \inf\{x : F^T(x) > 0\}$ and $s'_* = \inf\{x : F^S(x) > 0\}$ play the role of the right endpoints. Consequently, the results of Theorems~\ref{th:lim_corr} and~\ref{th:lim_coef} hold analogously for $u \to -\infty$ with respect to the left tail behaviour.

Future research may address these instabilities by exploring alternative coefficient selection strategies, such as leveraging robust estimates obtained at smaller $u$ values rather than direct estimation in extreme regimes.

\newpage
\section*{Data availability}
No external datasets were used in this study. All simulation data and results presented in this work can be generated and reproduced using the source code available at \url{https://github.com/LouB-N/Control-variates-MGF.git}.

\section*{Funding}
Part of Naveau's research work was supported by the Agence Nationale de la Recherche via: the SICIM and SHARE PEPR Maths-Vives project (France 2030 ANR-24-EXMA-0008), EXSTA grant (ANR-23-CE40-0009-01), VALPRED, PORC-EPIC, the PEPR TRACCS program (PC4 EXTENDING, ANR-22-EXTR-0005), and the PEPR IRIMONT (France 2030 ANR-22-EXIR-0003). He has also benefited from the Geolearning research chair.

\section*{Declaration of generative AI in the manuscript preparation process.}
During the preparation of this work, the authors used Gemini 3.5 to assist with language editing and code implementation. After using this tool, the authors reviewed and edited the content as needed and take full responsibility for the content of the published article.

\bibliographystyle{apalike}

\bibliography{cas-refs}

\makeatletter
\write\@auxout{\string\csxdef{mainlastpage}{\thepage}}
\makeatother
\clearpage
\appendix
\setcounter{page}{1}
\renewcommand{\thepage}{A\arabic{page}}
\def\lastpage{\applastpage}

\clearpage
\appendix
\setcounter{page}{1}
\renewcommand{\thepage}{A\arabic{page}}
\counterwithin{equation}{section}
\renewcommand{\theequation}{A.\arabic{equation}}

\section{Proofs}\label{sec:proofs}
\subsection{Asymptotic equivalent of $\E[\exp(u X)]$}
Let $X$ be a random variable in the max-domain of attraction of the Weibull distribution with EVI $\gamma < 0$ and with finite right endpoint $x^*=\text{sup}\{x:F_X(x)<1\}$ \citep{de_haan_extreme_2006}.\\ 

The cdf of $Y = x^* - X$ is regularly varying at 0 with index $-1/\gamma$, which is equivalent to the following:
\begin{align}
    &&\underset{t \to 0}{\lim} \frac{1-F_X(x^* - t x)}{1-F_X(x^* - t)} &= x^{-1/\gamma} \qquad \text{\cite[(1.2.3) from Theorem 1.2.1]{de_haan_extreme_2006}} \nonumber \\
    \Longleftrightarrow&& \underset{t \to 0}{\lim} \frac{1-\mathbb{P}(X \leq x^* - t x)}{1-\mathbb{P}(X \leq x^* - t)} &= x^{-1/\gamma}  \nonumber \\
    \Longleftrightarrow&& \underset{t \to 0}{\lim} \frac{1-\mathbb{P}(t x \leq x^* - X)}{1-\mathbb{P}(t \leq x^* - X)} &= x^{-1/\gamma}  \nonumber \\
    \Longleftrightarrow&& \underset{s \to +\infty}{\lim} \frac{1-\mathbb{P}\left(\frac{x}{s} \leq x^* - X\right)}{1-\mathbb{P}\left(\frac{1}{s} \leq x^* - X\right)} &= x^{-1/\gamma} \qquad \left(\text{with $t = \frac{1}{s}$} \right) \nonumber \\
    \Longleftrightarrow&& \underset{s \to +\infty}{\lim} \frac{\mathbb{P}\left(x^* - X \leq \frac{x}{s} \right)}{\mathbb{P}\left(x^* - X \leq \frac{1}{s}\right)} &= x^{-1/\gamma} \nonumber \\
    \underset{\bigstar}{\Longleftrightarrow}&& \underset{s \to +\infty}{\lim} \frac{\mathbb{P}\left(x^* - X \leq \frac{1}{s v} \right)}{\mathbb{P}\left(x^* - X \leq \frac{1}{s}\right)} &= v^{1/\gamma}  \qquad \left(\text{with $v = \frac{1}{x}$} \right) \nonumber \\
    \Longleftrightarrow&& \underset{s \to +\infty}{\lim} \frac{F_{x^* - X}\left(\frac{1}{s v} \right)}{F_{x^* - X}\left(\frac{1}{s}\right)} &= v^{1/\gamma} \nonumber \\
    \Longrightarrow&& \boxed{F_{x^* - X}\left(\frac{1}{v} \right) \underset{v \to +\infty}{\sim} v^{1/\gamma} L(v)} &\quad \text{where $L$ is a slowly varying function at $+\infty$.} \label{eq:reg_var}
\end{align}

Similarly, the survival function of $\frac{1}{x^* - X}$ is regularly varying at infinity with index $1/\gamma$ as follows from the previous relation:
\begin{align}
    &&\underset{t \to 0}{\lim} \frac{1-F_X(x^* - t x)}{1-F_X(x^* - t)} &= x^{-1/\gamma} \qquad \text{\cite[(1.2.3) from Theorem 1.2.1]{de_haan_extreme_2006}}  \nonumber \\
    \underset{\bigstar}{\Longleftrightarrow}&& \underset{s \to +\infty}{\lim} \frac{\mathbb{P}\left(x^* - X \leq \frac{1}{s v} \right)}{\mathbb{P}\left(x^* - X \leq \frac{1}{s}\right)} &= v^{1/\gamma} \nonumber \\
    \Longleftrightarrow&& \underset{s \to +\infty}{\lim} \frac{\mathbb{P}\left(\frac{1}{x^* - X} > s v \right)}{\mathbb{P}\left(\frac{1}{x^* - X} > s\right)} &= v^{1/\gamma}  \nonumber \\
    \Longleftrightarrow&& \underset{s \to +\infty}{\lim} \frac{1-F_{\frac{1}{x^* - X}}\left(s v\right)}{1-F_{\frac{1}{x^* - X}}\left(s\right)} &= v^{1/\gamma}  \nonumber \\
    \Longrightarrow&& \boxed{1-F_{\frac{1}{x^* - X}}\left(v\right) \underset{v \to +\infty}{\sim} v^{1/\gamma} L(v)} &\quad \text{where $L$ is a slowly varying function at $+\infty$.} \label{eq:reg_var_inv}
\end{align}

\newpage
From this, an asymptotic equivalent of the moment-generating function of $X$ can be found:
\begin{align}
    \E[\exp(u X)] &= \exp(u x^*) \E[\exp(-u Y)]&& \nonumber \\
    &= \exp(u x^*) \int_0^{+\infty}\exp(-u y)dF_Y(y)&& \nonumber \\
    &= \exp(u x^*) u \int_0^{+\infty}\exp(-u y)F_Y(y)dy \qquad &&\text{(by IPP)} \nonumber \\
    &= \exp(u x^*) \int_0^{+\infty}\exp(-t)F_Y\left(\frac{t}{u}\right)du \qquad &&\text{(with $t=u y$)} \nonumber \\
    &\underset{u \to +\infty}{\sim} \exp(u x^*) \int_0^{+\infty}\exp(-t)\left(\frac{u}{t}\right)^{\frac{1}{\gamma}} L\left(\frac{u}{t}\right) du \qquad &&\left(\text{by \eqref{eq:reg_var}, $F_Y\left(\frac{t}{u}\right)\underset{u \to +\infty}{\sim}\left(\frac{u}{t}\right)^{\frac{1}{\gamma}} L\left(\frac{u}{t}\right)$} \right) \nonumber \\
    &\underset{u \to +\infty}{\sim} \exp(u x^*) u^{\frac{1}{\gamma}} L\left(u\right) \int_0^{+\infty} \exp(-t) t^{-\frac{1}{\gamma}} du \qquad &&\left(\text{as L slowly varying, $L\left(\frac{u}{t}\right)\underset{u \to +\infty}{\sim}L\left(u\right)$}\right) \nonumber \\
    &\underset{u \to +\infty}{\sim} \exp(u x^*) u^{\frac{1}{\gamma}} L\left(u\right) \Gamma\left(1-\frac{1}{\gamma}\right).&& \label{eq:mgf_equivalent}
\end{align}

\vspace{10pt}
\subsection{Asymptotic equivalent of $\E\left[\exp\left(u\left(T+S\right)\right)\right]$}
Let $V^T=t^*-T$, $V^S=s^*-S$, and $u>0$.
\begin{align*}
    \E\left[\exp\left(u\left(T+S\right)\right)\right] &= \exp\left(u\left(t^*+s^*\right)\right) \E[\exp\left(- u\left(V^T+V^S\right)\right)]\\
    &= \exp\left(u\left(t^*+s^*\right)\right) \int_0^1 \mathbb{P}\left( \exp\left(- u\left(V^T+V^S\right)\right) >t \right)dt \\
    &= \exp\left(u\left(t^*+s^*\right)\right) \int_0^1 \mathbb{P}\left( V^T+V^S < \frac{-ln(t)}{u} \right)dt \\
    &= \exp\left(u\left(t^*+s^*\right)\right) \int_0^1 \mathbb{P}\left( \frac{1}{V^T+V^S} > \frac{-u}{ln(t)} \right)dt\\
    &= \exp\left(u\left(t^*+s^*\right)\right) \int_0^1 \mathbb{P}\left( \frac{1}{\frac{1}{Z^T}+\frac{1}{Z^S}} > \frac{-u}{ln(t)} \right)dt
\end{align*}
with $Z^T=\frac{1}{t^*-T}$ and $Z^S=\frac{1}{s^*-S}$.\\

By regular variations \eqref{eq:reg_var_inv}, $\mathbb{P}\left(Z^T > u\right) \underset{u \to +\infty}{\sim} u^{1 / \gamma^T} L_T(u) = u^{-\alpha_T} L_T(u)$, with $\alpha_T = -1 / \gamma^T$,\\ and $\mathbb{P}\left(Z^S > u\right) \underset{u \to +\infty}{\sim} u^{1 / \gamma^S} L_S(u) = u^{-\alpha_S} L_S(u)$, with $\alpha_S = -1 / \gamma^S$.\\

Assume that the joint tail probability satisfies a multivariate regular variation condition \citep{basrak2025foundations}, such that there exists a positive function $g$ going to 0 at infinity, and a measure $\nu(\cdot)$ defined on the unit sphere, for which \begin{align}\underset{u \to +\infty}{lim} \frac{\mathbb{P}\left( \left(\frac{1}{u \left(t^*-T\right)},\frac{1}{u \left(s^*-S\right)}\right) \in B \right)}{g(u)} = \nu(B) \ \text{ with } B = \left\{ (x,y)\in(0,\infty)^2; \frac{1}{x} + \frac{1}{y} < 1 \right\} .\end{align} \\

Assuming asymptotic dependence, if $\nu\left( (0, \infty)^2\right)>0$ then $\frac{1}{\frac{1}{Z^T}+\frac{1}{Z^S}}$ is regularly varying with 
\begin{align}
    &g(u) = \begin{cases}
        \mathbb{P}\left(Z^T>u\right) \text{ when $Z^T$ has lighter tail than $Z^S$};\\
        \mathbb{P}\left(Z^S>u\right) \text{ when $Z^S$ has lighter tail than $Z^T$}. \\
    \end{cases} \nonumber \\
    \Longrightarrow \quad &\mathbb{P}\left( \frac{1}{\frac{1}{Z^T}+\frac{1}{Z^S}} > \frac{-u}{ln(t)} \right) \underset{u \to +\infty}{\sim} \begin{cases}
        \nu(B) \; \mathbb{P}\left(Z^T>\frac{-u}{ln(t)}\right) \text{ when $Z^T$ has lighter tail than $Z^S$}; \\
        \nu(B) \; \mathbb{P}\left(Z^S>\frac{-u}{ln(t)}\right) \text{ when $Z^S$ has lighter tail than $Z^T$}. \\
    \end{cases} \nonumber \\
    \Longrightarrow \quad &\mathbb{P}\left( \frac{1}{\frac{1}{Z^T}+\frac{1}{Z^S}} > \frac{-u}{ln(t)} \right) \underset{u \to +\infty}{\sim} \nu(B) u^{-\alpha_*} \left(\frac{-1}{ln(t)}\right)^{-\alpha_*} \quad \text{ where } \alpha_*=max(\alpha_T,\alpha_S) \nonumber \\
    \Longrightarrow \quad & \E\left[\exp\left(u\left(T+S\right)\right)\right] \underset{u \to +\infty}{\sim} \exp\left(u\left(t^*+s^*\right)\right) \; \nu(B) \; u^{-\alpha_*} \; \left(\int_0^1 (-ln(t))^{\alpha_*}dt\right) \nonumber \\
    \Longrightarrow \quad & \E\left[\exp\left(u\left(T+S\right)\right)\right] \underset{u \to +\infty}{\sim} \exp\left(u\left(t^*+s^*\right)\right) \; \nu(B) \; u^{-\alpha_*} \; \Gamma\left(1+\alpha_*\right).\label{eq:mgf_equivalent_2}
\end{align}

In the case of asymptotic independence, i.e. if $\nu\left( (0, \infty)^2\right)=0$, then \begin{align}\E\left[\exp\left(u\left(T+S\right)\right)\right] \underset{u \to +\infty}{\sim} 0 . \label{eq:mgf_equivalent_2_ind}\end{align}

\vspace{10pt}
\subsection{Proof of Theorem~\ref{th:lim_corr}: limit of the correlation between $\tilde{T}(u)$ and $\tilde{S}(u)$} \label{sec:proof_lim_corr}
In the case of asymptotic dependence, joining \eqref{eq:mgf_equivalent} and \eqref{eq:mgf_equivalent_2} yields
\begin{align*}
    \corr\left(\tilde{T}(u),\tilde{S}(u)\right) &= \frac{\cov\left(\exp\left(u T\right), \exp\left(u S\right)\right)}{\sqrt{\var\left(\exp\left(u S\right)\right) \var\left(\exp\left(u T\right)\right)}}\\
    &= \frac{\E\left[\exp\left(u\left(S+T\right)\right)\right]-\E\left[\exp\left(u T\right)]\E[\exp\left(u S\right)\right]}{\sqrt{\left(\E\left[\exp\left(2 u S\right)\right]-\E\left[\exp\left(u S\right)\right]^2\right)\left(\E\left[\exp\left(2 u T\right)\right]-\E\left[\exp\left(u T\right)\right]^2\right)}} . \\
\end{align*}
Since the asymptotics are governed by exponential and power-law terms, the slowly varying functions $L_T$ and $L_S$ do not affect the leading-order behaviour and can be omitted at that level of approximation. Assuming $\lim_{u \to +\infty} L_T(u) = c_T > 0$ and $\lim_{u \to +\infty} L_S(u) = c_S > 0$.\\
Let $C_T = L_T(u) \Gamma\left(1-1 / \gamma^T\right) \underset{u \to +\infty}{\sim} c_T \Gamma\left(1-1 / \gamma^T\right)>0$,\\
and $C_S = L_S(u) \Gamma\left(1-1 / \gamma^S\right) \underset{u \to +\infty}{\sim} c_S \Gamma\left(1-1 / \gamma^S\right)>0$,\\
and $C_{TS} = \nu(B) \Gamma\left(1+\alpha_*\right) >0$.\\
\begin{align*}
    \corr\left(\tilde{T}(u),\tilde{S}(u)\right) &\underset{u \to +\infty}{\sim} \frac{\exp\left(u\left(t^*+s^*\right)\right) u^{-\alpha_*}C_{TS} - \exp\left(u t^*\right) u^{1 / \gamma^T} C_T  \exp\left(u s^*\right) u^{1 / \gamma^S} C_S}{ \sqrt{\left(\exp\left(2 u s^*\right) \left( (2 u)^{1 / \gamma^S} C_S - u^{2/\gamma^S}  C_S^2\right) \right) \; \left(\exp\left(2 u t^*\right) \left( (2 u)^{1 / \gamma^T} C_T - u^{2/\gamma^T}  C_T^2\right) \right)}} \\
    &\underset{u \to +\infty}{\sim} \frac{u^{-max\left(-1 / \gamma^T,-1 / \gamma^S\right)}C_{TS}}{u^{\frac{1}{2\gamma^S}+\frac{1}{2\gamma^T}}\sqrt{2^{1 / \gamma^S} C_S 2^{1 / \gamma^T} C_T}} \qquad \text{(considering only the higher-order terms)}\\
    &\underset{u \to +\infty}{\longrightarrow} \begin{cases}
        0 & \text{if } \gamma^T \neq \gamma^S;\\
        \frac{C_{TS}}{\sqrt{C_T C_S}}2^{-1 / (2\gamma^S) -1 / (2\gamma^T)} &  \text{if } \gamma^T = \gamma^S.
    \end{cases} \\
    &\underset{u \to +\infty}{\longrightarrow} \begin{cases}
        0 & \text{if } \gamma^T \neq \gamma^S;\\
        \frac{\nu(B)\Gamma\left(1+\alpha_*\right)}{\sqrt{c_T \Gamma\left(1-1 / \gamma^T\right) c_S \Gamma\left(1-1 / \gamma^S\right)}} 2^{-1 / (2\gamma^S) -1 / (2\gamma^T)} &  \text{if } \gamma^T = \gamma^S.
    \end{cases} \\
        &\underset{u \to +\infty}{\longrightarrow} \begin{cases}
        0 & \text{if } \gamma^T \neq \gamma^S;\\
        \frac{\nu(B)}{\sqrt{c_T c_S}} 2^{-1 / \gamma^T} &  \text{if } \gamma^T = \gamma^S.
    \end{cases}
\end{align*}

In the case of asymptotic independence, joining \eqref{eq:mgf_equivalent} and \eqref{eq:mgf_equivalent_2_ind} yields
\begin{align*}
    \corr\left(\tilde{T}(u),\tilde{S}(u)\right) &= \underset{u \to +\infty}{\sim} \frac{- \exp\left(u\left(t^*+s^*\right)\right)  u^{1 / \gamma^T} C_T  u^{1 / \gamma^S} C_S}{ \sqrt{\left(\exp\left(2 u s^*\right) \left( (2 u)^{1 / \gamma^S} C_S - u^{2/\gamma^S}  C_S^2\right) \right) \; \left(\exp\left(2 u t^*\right) \left( (2 u)^{1 / \gamma^T} C_T - u^{2/\gamma^T}  C_T^2\right) \right)}}\\
    &\underset{u \to +\infty}{\sim} - u^{1/(2\gamma^T) + 1/(2\gamma^S)} \frac{C_T C_S}{\sqrt{2^{1 / \gamma^S + 1 / \gamma^T} C_S C_T}} \qquad \text{(considering only the higher-order terms)}\\
    &\underset{u \to +\infty}{\longrightarrow} 0 .
\end{align*}

\vspace{10pt}
\subsection{Proof of Theorem~\ref{th:lim_coef}: limit of $w_{opt}(u)$} \label{sec:proof_lim_coef}
In the case of asymptotic dependence, joining \eqref{eq:mgf_equivalent} and \eqref{eq:mgf_equivalent_2} yields
\begin{align*}
    w_{opt}(u) &= \frac{\cov\left(\exp\left(u T\right), \exp\left(u S\right)\right)}{\var\left(\exp\left(u S\right)\right)} = \frac{\E\left[\exp\left(u\left(S+T\right)\right)\right]-\E\left[\exp\left(u T\right)\right]\E\left[\exp\left(u S\right)\right]}{\E\left[\exp\left(2 u S\right)\right]-\E\left[\exp\left(u S\right)\right]^2} \\
    &\underset{u \to +\infty}{\sim} \exp\left(u\left(t^*-s^*\right)\right) \frac{u^{-\alpha_*}\nu(B) \Gamma\left(1+\alpha_*\right) - u^{1 / \gamma^T} L_T(u) \Gamma\left(1-1 / \gamma^T \right)  u^{1 / \gamma^S} L_S(u) \Gamma\left(1-1 / \gamma^S \right)}{(2 u)^{1 / \gamma^S} L_S(2 u) \Gamma\left(1-1 / \gamma^S \right) - u^{2/\gamma^S} L_S^2(u) \Gamma\left(1-1 / \gamma^S \right)^2} .
\end{align*}
Since the asymptotics are governed by exponential and power-law terms, the slowly varying functions $L_T$ and $L_S$ do not affect the leading-order behaviour and can be omitted at that level of approximation.Assuming $\lim_{u \to +\infty} L_T(u) = c_T > 0$ and $\lim_{u \to +\infty} L_S(u) = c_S > 0$.\\
Let $C_T = L_T(u) \Gamma\left(1-1 / \gamma^T\right) \underset{u \to +\infty}{\sim} c_T \Gamma\left(1-1 / \gamma^T\right)>0$,\\
and $C_S = L_S(u) \Gamma\left(1-1 / \gamma^S\right) \underset{u \to +\infty}{\sim} c_S \Gamma\left(1-1 / \gamma^S\right)>0$,\\
and $C_{TS} = \nu(B) \Gamma\left(1+\alpha_*\right) >0$.\\
\begin{align*}
    w_{opt}(u) &\underset{u \to +\infty}{\sim} \exp\left(u\left(t^*-s^*\right)\right) \frac{u^{-\alpha_*}C_{TS} - u^{1 / \gamma^T} C_T  u^{1 / \gamma^S} C_S}{(2 u)^{1 / \gamma^S} C_S - u^{2/\gamma^S}  C_S^2}\\
    &\underset{u \to +\infty}{\sim} \exp\left(u\left(t^*-s^*\right)\right) \frac{u^{-\alpha_* -1 / \gamma^S}C_{TS}}{2^{1 / \gamma^S} C_S} \qquad \text{(considering only the higher-order terms)}\\
    &\underset{u \to +\infty}{\longrightarrow} \begin{cases}
        + \infty & \text{if } t^*>s^*;\\
        \begin{cases}
            \frac{C_{TS}}{2^{1 / \gamma^S} C_S} = \frac{\nu(B)}{2^{1 / \gamma^S} c_S} & \text{if } \gamma^T \leq \gamma^S;\\
            0 & \text{if } \gamma^T > \gamma^S;\\
        \end{cases} & \text{and } t^*=s^*; \\
        0 &  \text{if } t^*<s^*;
    \end{cases}
\end{align*}
as the exponent $-\alpha_* - 1 / \gamma^S = - \max\left(-1 / \gamma^T, -1 / \gamma^S\right) - 1 / \gamma^S = \begin{cases}
    0 & \text{if } \gamma^T \leq \gamma^S; \\
    1 / \gamma^T - 1 / \gamma^S <0 & \text{if } \gamma^T > \gamma^S.\\
\end{cases}$\\

Under asymptotic dependence, i.e. $\nu\left( (0, \infty)^2\right)>0$, the asymptotic behaviour of $w_{opt}(u)$ reflects the tail behaviours of $\tilde{T}(u)$ and $\tilde{S}(u)$. When $t^* \neq s^*$, the variable with the larger endpoint dominates: $t^* > s^*$ drives $w_{opt}(u) \to +\infty$, while $t^* < s^*$ drives $w_{opt}(u) \to 0$. Equal endpoints are thus necessary for stabilization, though not sufficient. In this case, relative tail decay is governed by the EVIs. If $\gamma^T > \gamma^S$, $\tilde{S}(u)$ decays faster, driving $w_{opt}(u) \to 0$. Only when $\gamma^T \leq \gamma^S$, the coefficient $w_{opt}(u)$ stabilizes at a finite non-zero constant. However, a stable optimal coefficient $w_{opt}(u)$ does not guarantee variance reduction, as seen with the limit of the correlation, that is non-zero if the EVIs are equal.\\

In the case of asymptotic independence, joining \eqref{eq:mgf_equivalent} and \eqref{eq:mgf_equivalent_2_ind} yields
\begin{align*}
    w_{opt}(u) &\underset{u \to +\infty}{\sim} -\exp\left(u\left(t^*-s^*\right)\right) \frac{u^{1 / \gamma^T} L_T(u) \Gamma\left(1-1 / \gamma^T \right) u^{1 / \gamma^S} L_S(u) \Gamma\left(1-1 / \gamma^S \right)}{(2 u)^{1 / \gamma^S} L_S(2 u) \Gamma\left(1-1 / \gamma^S \right) - u^{2/\gamma^S} L_S^2(u) \Gamma\left(1-1 / \gamma^S \right)^2} \\
     &\underset{u \to +\infty}{\sim} -\exp\left(u\left(t^*-s^*\right)\right) u^{1 / \gamma^T} \frac{\Gamma\left(1-1 / \gamma^T \right)}{2^{1 / \gamma^S}} \qquad \text{(considering only the higher-order terms)}\\
     &\underset{u \to +\infty}{\longrightarrow} \begin{cases}
            - \infty & \text{if } t^* > s^*;\\
            0 & \text{if } t^* \leq s^*. \\
        \end{cases}
\end{align*}

Under asymptotic independence, i.e. $\nu\left( (0, \infty)^2\right)=0$, the joint tail co-movement vanishes as $u \to +\infty$. In this regime, $w_{opt}(u)$ cannot settle at a finite non-zero constant: if $t^* > s^*$, the scale discrepancy forces $w_{opt}(u) \to -\infty$, whereas if $t^* \leq s^*$, $w_{opt}(u)$ vanishes to $0$.

\end{document}